\documentclass{IEEEtran}

\usepackage{amsmath}
\usepackage{amssymb}
\usepackage{cite}
\usepackage{graphicx}
\usepackage{verbatim}
\usepackage{subfigure}
\usepackage{color}

\newcommand{\mb}[1]{{\boldsymbol #1}}

\renewcommand{\b}{\mb{b}}

\newcommand{\e}{\mb{e}}
\newcommand{\f}{\mb{f}}
\newcommand{\g}{\mb{g}}
\newcommand{\h}{\mb{h}}

\renewcommand{\r}{\mb{r}}
\renewcommand{\u}{\mb{u}}
\renewcommand{\v}{\mb{v}}
\newcommand{\x}{\mb{x}}
\newcommand{\w}{\mb{w}}
\newcommand{\y}{\mb{y}}

\newcommand{\A}{\mb{A}}
\newcommand{\B}{\mb{B}}
\newcommand{\CC}{{\mathcal C}}
\newcommand{\D}{\mb{D}}
\newcommand{\F}{\mb{F}}
\newcommand{\FF}{\mathcal F}
\renewcommand{\H}{\mb{H}}
\newcommand{\I}{\mb{I}}
\newcommand{\J}{\mb{J}}
\newcommand{\M}{\mb{M}}

\renewcommand{\SS}{{\mathcal S}}
\newcommand{\U}{\mb{U}}

\renewcommand{\P}{\mb{P}}
\newcommand{\Q}{\mb{Q}}
\newcommand{\TT}{{\mathcal T}}
\newcommand{\W}{\mb{W}}
\newcommand{\XX}{{\mathcal X}}

\newcommand{\bD}{{\boldsymbol \Delta}}

\newcommand{\zero}{\mb{0}}
\newcommand{\eps}{\varepsilon}
\newcommand{\Bxz}{{{\mathcal B}_\eps(\x_0)}}
\newcommand{\Baz}{{{\mathcal B}_\eps(\alf_0)}}
\newcommand{\st}{\text{s.t. }}
\newcommand{\Hal}{\H_\alf}
\newcommand{\Halz}{\H_{\alf_0}}

\newcommand{\UUJUU}{{\U\U^T\J\U\U^T}}

\newcommand{\beq}{\begin{equation}}
\newcommand{\eeq}{\end{equation}}

\newcommand{\hx}{{\widehat{\x}}}
\newcommand{\xu}{\x_\mathrm{u}}
\newcommand{\alf}{{\boldsymbol \alpha}}
\newcommand{\half}{{\widehat{\alf}}}
\newcommand{\alfor}{\half_{\mathrm{oracle}}}
\newcommand{\alfds}{\half_{\mathrm{DS}}}
\newcommand{\alfgds}{\half_{\mathrm{GDS}}}
\newcommand{\alfbp}{\half_{\mathrm{BP}}}

\newcommand{\RR}{{\mathbb R}}

\newcommand{\MSE}{{\mathrm{MSE}}}

\newcommand{\E}[1]{E \! \left\{ #1 \right\}}
\newcommand{\pd}[2]{\frac{\partial #1}{\partial #2}}

\newcommand{\Nu}[1]{{{\mathcal N}\! \left( #1 \right) }}
\newcommand{\Ra}[1]{{{\mathcal R} ( #1 ) }}

\newcommand{\pinv}{\dagger}

\newtheorem{theorem}{Theorem}

\DeclareMathOperator{\Cov}{Cov}
\DeclareMathOperator{\Tr}{Tr}

\DeclareMathOperator{\spark}{spark}
\DeclareMathOperator{\spn}{span}
\DeclareMathOperator{\supp}{supp}

\begin{document}

\title{The~Cram\'er--Rao~Bound for~Sparse~Estimation}
\author{Zvika~Ben-Haim,~\IEEEmembership{Student~Member,~IEEE,}
and~Yonina~C.~Eldar,~\IEEEmembership{Senior~Member,~IEEE}%
\thanks{Department of Electrical Engineering, Technion---Israel Institute of Technology, Haifa 32000, Israel. Phone: +972-4-8294700, fax: +972-4-8295757, E-mail: \{zvikabh@tx, yonina@ee\}.technion.ac.il. This work was supported in part by the Israel Science Foundation under Grant no. 1081/07 and by the European Commission in the framework of the FP7 Network of Excellence in Wireless COMmunications NEWCOM++ (contract no. 216715).}}
\maketitle

\begin{abstract}
The goal of this paper is to characterize the best achievable performance for the problem of estimating an unknown parameter having a sparse representation. Specifically, we consider the setting in which a sparsely representable deterministic parameter vector is to be estimated from measurements corrupted by Gaussian noise, and derive a lower bound on the mean-squared error (MSE) achievable in this setting. To this end, an appropriate definition of bias in the sparse setting is developed, and the constrained Cram\'er--Rao bound (CRB) is obtained. This bound is shown to equal the CRB of an estimator with knowledge of the support set, for almost all feasible parameter values. Consequently, in the unbiased case, our bound is identical to the MSE of the oracle estimator. Combined with the fact that the CRB is achieved at high signal-to-noise ratios by the maximum likelihood technique, our result provides a new interpretation for the common practice of using the oracle estimator as a gold standard against which practical approaches are compared.
\end{abstract}

{\it EDICS Topics:} SSP-PARE, SSP-PERF.


\section{Introduction}

The problem of estimating a sparse unknown parameter vector from noisy measurements has been analyzed intensively in the past few years \cite{tropp06, donoho06, candes06, candes07}, and has already given rise to numerous successful signal processing algorithms \cite{elad06, dabov07, dabov08, protter09, elad05}. In this paper, we consider the setting in which noisy measurements of a deterministic vector $\x_0$ are available. It is assumed that $\x_0$ has a sparse representation $\x_0 = \D\alf_0$, where $\D$ is a given dictionary and most of the entries of $\alf_0$ equal zero. Thus, only a small number of ``atoms,'' or columns of $\D$, are required to represent $\x_0$. The challenges confronting an estimation technique are to recover either $\x_0$ itself or its sparse representation $\alf_0$. Several practical approaches turn out to be surprisingly successful in this task. Such approaches include the Dantzig selector (DS) \cite{candes07} and basis pursuit denoising (BPDN), which is also referred to as the Lasso \cite{chen98, tropp06, donoho06}.

A standard measure of estimator performance is the mean-squared error (MSE). Several recent papers analyzed the MSE obtained by methods such as the DS and BPDN \cite{candes07, ben-haim09c}. To determine the quality of estimation approaches, it is of interest to compare their achievements with theoretical performance limits: if existing methods approach the performance bound, then they are nearly optimal and further improvements in the current setting are impossible. This motivates the development of lower bounds on the MSE of estimators in the sparse setting.

Since the parameter to be estimated is deterministic, the MSE is in general a function of the parameter value. While there are lower bounds on the worst-case achievable MSE among all possible parameter values \cite[\S7.4]{candes06b}, the actual performance for a specific value, or even for most values, might be substantially lower. Our goal is therefore to characterize the minimum MSE obtainable for each particular parameter vector. A standard method of achieving this objective is the Cram\'er--Rao bound (CRB) \cite{kay93, shao03}.

The fact that $\x_0$ has a sparse representation is of central importance for estimator design. Indeed, many sparse estimation settings are underdetermined, meaning that without the assumption of sparsity, it is impossible to identify the correct parameter from its measurements, even without noise. In this paper, we treat the sparsity assumption as a deterministic prior constraint on the parameter. Specifically, we assume that $\x_0 \in \SS$, where $\SS$ is the set of all parameter vectors which can be represented by no more than $s$ atoms, for a given integer $s$.

Our results are inspired by the well-studied theory of the constrained CRB \cite{gorman90, marzetta93, StoicaNg98, ben-haim09}. This theory is based on the assumption that the constraint set can be defined using the system of equations $\f(\x)=\zero$, $\g(\x)\le\zero$, where $\f$ and $\g$ are continuously differentiable functions. The resulting bound depends on the derivatives of the function $\f$. However, sparsity constraints cannot be written in this form. This necessitates the development of a bound suitable for non-smooth constraint sets \cite{ben-haim09d}. In obtaining this modified bound, we also provide new insight into the meaning of the general constrained CRB\@. In particular, we show that the fact that the constrained CRB is lower than the unconstrained bound results from an expansion of the class of estimators under consideration.

With the aforementioned theoretical tools at hand, we obtain lower bounds on the MSE in a variety of sparse estimation problems. Our bound limits the MSE achievable by any estimator having a pre-specified bias function, for each parameter value. Particular emphasis is given to the unbiased case; the reason for this preference is twofold: First, when the signal-to-noise ratio (SNR) is high, biased estimation is suboptimal. Second, for high SNR values, the unbiased CRB is achieved by the maximum likelihood (ML) estimator.

While the obtained bounds differ depending on the exact problem definition, in general terms and for unbiased estimation the bounds can be described as follows. For parameters having maximal support, i.e., parameters whose representation requires the maximum allowed number $s$ of atoms, the lower bound equals the MSE of the ``oracle estimator'' which knows the locations (but not the values) of the nonzero representation elements. On the other hand, for parameters which do not have maximal support (a set which has Lebesgue measure zero in $\SS$), our lower bound is identical to the CRB for an unconstrained problem, which is substantially higher than the oracle MSE\@.

The correspondence between the CRB and the MSE of the oracle estimator (for all but a zero-measure subset of the feasible parameter set $\SS$) is of practical interest since, unlike the oracle estimator, the CRB is achieved by the ML estimator at high SNR\@. Our bound can thus be viewed as an alternative justification for the common use of the oracle estimator as a baseline against which practical algorithms are compared. This gives further merit to recent results, which demonstrate that BPDN and the DS both achieve near-oracle performance \cite{candes07, ben-haim09c}. However, the existence of parameters for which the bound is much higher indicates that oracular performance cannot be attained for \emph{all} parameter values, at least using unbiased techniques. Indeed, as we will show, in many sparse estimation scenarios, one cannot construct \emph{any} estimator which is unbiased for all sparsely representable parameters.

Our contribution is related to, but distinct from, the work of Babadi et al.\ \cite{babadi09}, in which the CRB of the oracle estimator was derived (and shown to equal the aforementioned oracle MSE). Our goal in this work is to obtain a lower bound on the performance of estimators which are not endowed with oracular knowledge; consequently, as explained above, for some parameter values the obtained CRB will be higher than the oracle MSE\@. It was further shown in \cite{babadi09} that when the measurements consist of Gaussian random mixtures of the parameter vector, there exists an estimator which achieves the oracle CRB at high SNR; this is shown to hold on average over realizations of the measurement mixtures. The present contribution strengthens this result by showing that for any given (deterministic) well-behaved measurement setup, there exists a technique (namely, the ML estimator) achieving the CRB at high SNR\@. Thus, convergence to the CRB is guaranteed for all measurement settings, and not merely when averaging over an ensemble of such settings.

The rest of this paper is organized as follows. In Section~\ref{se:sparse backgnd}, we review the sparse setting as a constrained estimation problem. Section~\ref{se:crb} defines a generalization of sparsity constraints, which we refer to as locally balanced constraint sets; the CRB is then derived in this general setting. In Section~\ref{se:sparse bounds}, our general results are applied back to some specific sparse estimation problems. In Section~\ref{se:numer}, the CRB is compared to the empirical performance of estimators of sparse vectors. Our conclusions are summarized in Section~\ref{se:discuss}.

Throughout the paper, boldface lowercase letters $\v$ denote vectors while boldface uppercase letters $\M$ denote matrices. Given a vector function $\f: \RR^n \rightarrow \RR^k$, we denote by $\partial \f / \partial \x$ the $k \times n$ matrix whose $ij$th element is $\partial f_i / \partial x_j$. The support of a vector, denoted $\supp(\v)$, is the set of indices of the nonzero entries in $\v$. The Euclidean norm of a vector $\v$ is denoted $\|\v\|_2$, and the number of nonzero entries in $\v$ is $\|\v\|_0$. Finally, the symbols $\Ra{\M}$, $\Nu{\M}$, and $\M^\pinv$ refer, respectively, to the column space, null space, and Moore--Penrose pseudoinverse of the matrix $\M$.

\section{Sparse Estimation Problems}
\label{se:sparse backgnd}

In this section, we describe several estimation problems whose common theme is that the unknown parameter has a sparse representation with respect to a known dictionary. We then review some standard techniques used to recover the unknown parameter in these problems. In Section~\ref{se:numer} we will compare these methods with the performance bounds we develop.

\subsection{The Sparse Setting}
\label{ss:sparse setting}

Suppose we observe a measurement vector $\y \in \RR^m$, given by
\beq \label{eq:y=Ax+w}
\y = \A\x_0 + \w
\eeq
where $\x_0 \in \RR^n$ is an unknown deterministic signal, $\w$ is independent, identically distributed (IID) Gaussian noise with zero mean and variance $\sigma^2$, and $\A$ is a known $m \times n$ matrix. We assume the prior knowledge that there exists a sparse representation of $\x_0$, or, more precisely, that
\beq \label{eq:def S}
\x_0 \in \SS \triangleq \left\{ \x \in \RR^n : \x = \D \alf, \|\alf\|_0 \le s \right\}.
\eeq
In other words, the set $\SS$ describes signals $\x$ which can be formed from a linear combination of no more than $s$ columns, or atoms, from $\D$. The dictionary $\D$ is an $n \times p$ matrix with $n \le p$, and we assume that $s < p$, so that only a subset of the atoms in $\D$ can be used to represent any signal in $\SS$. We further assume that $\D$ and $s$ are known.

Quite a few important signal recovery applications can be formulated using the setting described above. For example, if $\A=\I$, then $\y$ consists of noisy observations of $\x_0$, and recovering $\x_0$ is a denoising problem \cite{elad06, dabov07}. If $\A$ corresponds to a blurring kernel, we obtain a deblurring problem \cite{dabov08}. In both cases, the matrix $\A$ is square and invertible. Interpolation and inpainting can likewise be formulated as \eqref{eq:y=Ax+w}, but in those cases $\A$ is an underdetermined matrix, i.e., we have $m<n$ \cite{elad05}. For all of these estimation scenarios, our goal is to obtain an estimate $\hx$ whose MSE is as low as possible, where the MSE is defined as
\beq \label{eq:def MSE}
\MSE \triangleq \E{ \| \hx - \x_0 \|_2^2 }.
\eeq
Note that $\x_0$ is deterministic, so that the expectation in \eqref{eq:def MSE} (and throughout the paper) is taken over the noise $\w$ but not over $\x_0$. Thus, the MSE is in general a function of $\x_0$.

In the above settings, the goal is to estimate the unknown signal $\x_0$. However, it may also be of interest to recover the coefficient vector $\alf_0$ for which $\x_0 = \D\alf_0$, e.g., for the purpose of model selection \cite{tropp06, candes07}. In this case, the goal is to construct an estimator $\half$ whose MSE $E\{ \|\half-\alf_0\|_2^2 \}$ is as low as possible. Unless $\D$ is unitary, estimating $\alf_0$ is not equivalent to estimating $\x_0$. Note, however, that when estimating $\alf_0$, the matrices $\A$ and $\D$ can be combined to obtain the equivalent problem
\beq \label{eq:y=H alf + w}
\y = \H \alf_0 + \w
\eeq
where $\H \triangleq \A\D$ is an $m \times p$ matrix and
\beq \label{eq:def T}
\alf_0 \in \TT = \{ \alf \in \RR^p : \|\alf\|_0 \le s \} .
\eeq
Therefore, this problem can also be seen as a special case of \eqref{eq:y=Ax+w} and \eqref{eq:def S}. Nevertheless, it will occasionally be convenient to refer specifically to the problem of estimating $\alf_0$ from \eqref{eq:y=H alf + w}.

Signal estimation problems differ in the properties of the dictionary $\D$ and measurement matrix $\A$. In particular, problems of a very different nature arise depending on whether the dictionary is a basis or an overcomplete frame. For example, many approaches to denoising yield simple shrinkage techniques when $\D$ is a basis, but deteriorate to NP-hard optimization problems when $\D$ is overcomplete \cite{natarajan95}.

A final technical comment is in order. If the matrix $\H$ in \eqref{eq:y=H alf + w} does not have full column rank, then there may exist different feasible parameters $\alf_1$ and $\alf_2$ such that $\H\alf_1 = \H\alf_2$. In this case, the probability distribution of $\y$ will be identical for these two parameter vectors, and the estimation problem is said to be unidentifiable \cite[\S1.5.2]{lehmann98}. A necessary and sufficient condition for identifiability is
\beq \label{eq:spark req}
\spark(\H) > 2s
\eeq
where $\spark(\H)$ is defined as the smallest integer $k$ such that there exist $k$ linearly dependent columns in $\H$ \cite{donoho03}. We will adopt the assumption \eqref{eq:spark req} throughout the paper. Similarly, in the problem \eqref{eq:y=Ax+w} we will assume that
\beq \label{eq:spark req D}
\spark(\D) > 2s.
\eeq

\subsection{Estimation Techniques}
\label{ss:est techniques}

We now review some standard estimators for the sparse problems described above. These techniques are usually viewed as methods for obtaining an estimate $\half$ of the vector $\alf_0$ in \eqref{eq:y=H alf + w}, and we will adopt this perspective in the current section. One way to estimate $\x_0$ in the more general problem \eqref{eq:y=Ax+w} is to first estimate $\alf_0$ with the methods described below and then use the formula $\hx = \D\half$.

A widely-used estimation technique is the ML approach, which provides an estimate of $\alf_0$ by solving
\beq \label{eq:ml}
\min_\alf  \|\y - \H\alf\|_2^2 \quad \st \|\alf\|_0 \le s.
\eeq
Unfortunately, \eqref{eq:ml} is a nonconvex optimization problem and solving it is NP-hard \cite{natarajan95}, meaning that an efficient algorithm providing the ML estimator is unlikely to exist. In fact, to the best of our knowledge, the most efficient method for solving \eqref{eq:ml} for general $\H$ is to enumerate the $\binom{p}{s}$ possible $s$-element support sets of $\alf$ and choose the one for which $\|\y - \H\alf\|_2^2$ is minimal. This is clearly an impractical strategy for reasonable values of $p$ and $s$. Consequently, several efficient alternatives have been proposed for estimating $\alf_0$. One of these is the $\ell_1$-penalty version of BPDN \cite{tropp06}, which is defined as a solution $\alfbp$ to the quadratic program
\beq \label{eq:bpdn}
\min_\alf \tfrac{1}{2} \|\y - \H\alf\|_2^2 + \gamma \|\alf\|_1
\eeq
with some regularization parameter $\gamma$. More recently, the DS was proposed \cite{candes07}; this approach estimates $\alf_0$ as a solution $\alfds$ to
\beq \label{eq:ds}
\min_\alf \|\alf\|_1 \quad \st \|\H^T (\y - \H\alf) \|_\infty \le \tau
\eeq
where $\tau$ is again a user-selected parameter. A modification of the DS, known as the Gauss--Dantzig selector (GDS) \cite{candes07}, is to use $\alfds$ only to estimate the support of $\alf_0$. In this approach, one solves \eqref{eq:ds} and determines the support set of $\alfds$. The GDS estimate is then obtained as
\beq \label{eq:gds}
\alfgds =
\begin{cases}
\H_{\alfds}^\pinv \y    & \text{on the support set of $\alfds$} \cr
\zero                   & \text{elsewhere}
\end{cases}
\eeq
where $\H_{\alfds}$ consists of the columns of $\H$ corresponding to the support of $\alfds$.

Previous research on the performance of these estimators has primarily examined their worst-case MSE among all possible values of $\alf_0 \in \TT$. Specifically, it has been shown \cite{candes07} that, under suitable conditions on $\H$, $s$, and $\tau$, the DS of \eqref{eq:ds} satisfies
\beq \label{eq:ds wc bound}
\|\alf_0 - \alfds\|_2^2 \le C s \sigma^2 \log p \quad \text{with high probability}
\eeq
for some constant $C$. It follows that the MSE of the DS is also no greater than a constant times $s \sigma^2 \log p$ for all $\alf_0 \in \TT$ \cite{candes06b}. An identical property was also demonstrated for BPDN \eqref{eq:bpdn} with an appropriate choice of $\gamma$ \cite{ben-haim09c}. Conversely, it is known that the worst-case error of \emph{any} estimator is at least a constant times $s \sigma^2 \log p$ \cite[\S7.4]{candes06b}. Thus, both BPDN and the DS are optimal, up to a constant, in terms of worst-case error. Nevertheless, the MSE of these approaches for specific values of $\alf_0$, even for a vast majority of such values, might be much lower. Our goal differs from this line of work in that we characterize the \emph{pointwise} performance of an estimator, i.e., the MSE for specific values of $\alf_0$.

Another baseline with which practical techniques are often compared is the oracle estimator, given by
\beq\label{eq:def xo}
\alfor =
\begin{cases}
\Halz^\pinv \b & \text{on the set $\supp(\alf_0)$} \\
\zero         & \text{elsewhere}
\end{cases}
\eeq
where $\Halz$ is the submatrix constructed from the columns of $\H$ corresponding to the nonzero entries of $\alf_0$. In other words, $\alfor$ is the least-squares (LS) solution among vectors whose support coincides with $\supp(\alf_0)$, which is assumed to have been provided by an ``oracle.'' Of course, in practice the support of $\alf_0$ is unknown, so that $\alfor$ cannot actually be implemented. Nevertheless, one often compares the performance of true estimators with $\alfor$, whose MSE is given by \cite{candes07}
\beq \label{eq:oracle mse}
\sigma^2 \Tr((\Halz^T \Halz)^{-1}).
\eeq

Is \eqref{eq:oracle mse} a bound on estimation MSE\@? While $\alfor$ is a reasonable technique to adopt if $\supp(\alf_0)$ is known, this does not imply that \eqref{eq:oracle mse} is a lower bound on the performance of practical estimators. Indeed, as will be demonstrated in Section~\ref{se:numer}, when the SNR is low, both BPDN and the DS outperform $\alfor$, thanks to the use of shrinkage in these estimators. Furthermore, if $\supp(\alf_0)$ is known, then there exist biased techniques which are better than $\alfor$ for \emph{all} values of $\alf_0$ \cite{ben-haim06}. Thus, $\alfor$ is neither achievable in practice, nor optimal in terms of MSE\@. As we will see, one can indeed interpret \eqref{eq:oracle mse} as a lower bound on the achievable MSE, but such a result requires a certain restriction of the class of estimators under consideration.

\section{The Constrained Cram\'er--Rao Bound}
\label{se:crb}

A common technique for determining the achievable performance in a given estimation problem is to calculate the CRB, which is a lower bound on the MSE of estimators having a given bias \cite{kay93}. In this paper, we are interested in calculating the CRB when it is known that the parameter $\x$ satisfies sparsity constraints such as those of the sets $\SS$ of \eqref{eq:def S} and $\TT$ of \eqref{eq:def T}.

The CRB for constrained parameter sets has been studied extensively in the past \cite{gorman90, marzetta93, StoicaNg98, ben-haim09}. However, in prior work derivation of the CRB assumed that the constraint set is given by
\beq \label{eq:constr}
\XX = \{ \x \in \RR^n: \f(\x)=\zero, \ \g(\x)\le\zero \}
\eeq
where $\f(\x)$ and $\g(\x)$ are continuously differentiable functions. We will refer to such $\XX$ as continuously differentiable sets. As shown in prior work \cite{gorman90}, the resulting bound depends on the derivatives of the function $\f$. Yet in some cases, including the sparse estimation scenarios discussed in Section~\ref{se:sparse backgnd}, the constraint set cannot be written in the form \eqref{eq:constr}, and the aforementioned results are therefore inapplicable. Our goal in the current section is to close this gap by extending the constrained CRB to constraint sets $\XX$ encompassing the sparse estimation scenario.

We begin this section with a general discussion of the CRB and the class of estimators to which it applies. This will lead us to interpret the constrained CRB as a bound on estimators having an incompletely specified bias gradient. This interpretation will facilitate the application of the existing constrained CRB to the present context.

\begin{figure*}  
\centerline{\includegraphics{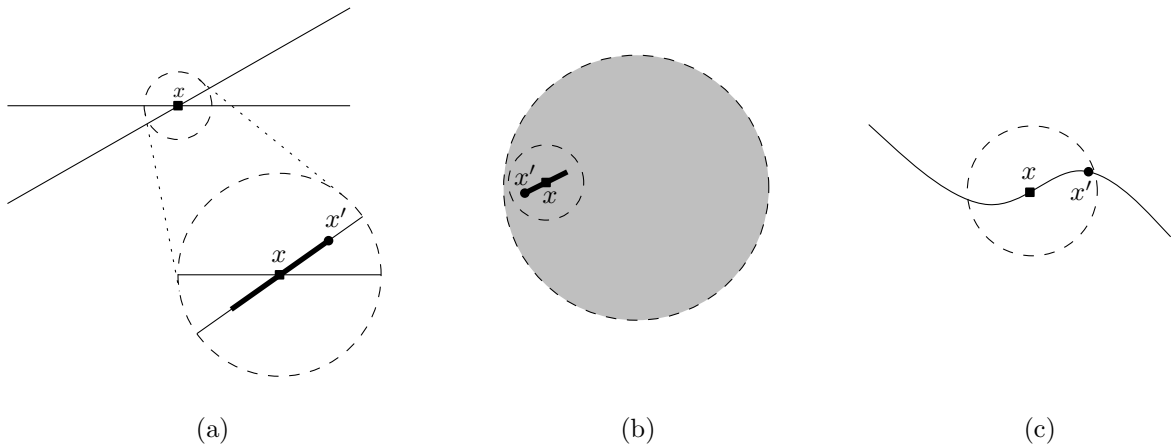}}
\caption{In a locally balanced set such as a union of subspaces (a) and an open ball (b), each point is locally defined by a set of feasible directions along which an infinitesimal movement does not violate the constraints. The curve (c) is not characterized in this way and thus is not locally balanced.}
\label{fi:loc bal}
\end{figure*}

\subsection{Bias Requirements in the Constrained CRB}
\label{ss:bias req}

In previous settings for which the constrained CRB was derived, it was noted that the resulting bound is typically lower than the unconstrained version \cite[Remark~4]{gorman90}. At first glance, one would attribute the reduction in the value of the CRB to the fact that the constraints add information about the unknown parameter, which can then improve estimation performance. On the other hand, the CRB separately characterizes the achievable performance for each value of the unknown parameter $\x_0$. Thus, the CRB at $\x_0$ applies even to estimators designed specifically to perform well at $\x_0$. Such estimators surely cannot achieve further gain in performance if it is known that $\x_0 \in \XX$. Why, then, is the constrained CRB lower than the unconstrained bound? The answer to this apparent paradox involves a careful definition of the class of estimators to which the bound applies.

To obtain a meaningful bound, one must exclude some estimators from consideration. Unless this is done, the bound will be tarnished by estimators of the type $\hx = \xu$, for some constant $\xu$, which achieve an MSE of $0$ at the specific point $\x = \xu$. It is standard practice to circumvent this difficulty by restricting attention to estimators having a particular bias $\b(\x) \triangleq \E{\hx} - \x$. In particular, it is common to examine unbiased estimators, for which $\b(\x) = \zero$.

However, in some settings, it is impossible to construct estimators which are unbiased for all $\x \in \RR^n$. For example, suppose we are to estimate the coefficients $\alf_0$ of an overcomplete dictionary based on the measurements given by \eqref{eq:y=H alf + w}. Since the dictionary is overcomplete, its nullspace is nontrivial; furthermore, each coefficient vector in the nullspace yields an identical distribution of the measurements, so that an estimator can be unbiased for one of these vectors at most.

The question is whether it is possible to construct estimators which are unbiased for some, but not all, values of $\x$. One possible approach is to seek estimators which are unbiased for all $\x \in \XX$. However, as we will see later in this section, even this requirement can be too strict: in some cases it is impossible to construct estimators which are unbiased for all $\x \in \XX$. More generally, the CRB is a \emph{local} bound, meaning that it determines the achievable performance at a particular value of $\x$ based on the statistics at $\x$ and at nearby values. Thus, it is irrelevant to introduce requirements on estimation performance for parameters which are distant from the value $\x$ of interest.

Since we seek a locally unbiased estimator, one possibility is to require unbiasedness at a single point, say $\xu$. As it turns out, it is always possible to construct such a technique: this is again $\hx = \xu$, which is unbiased at $\xu$ but nowhere else. To avoid this loophole, one can require an estimator to be unbiased in the neighborhood
\beq
\Bxz = \left\{ \x \in \RR^m : \|\x-\x_0\|_2 < \eps \right\}
\eeq
of $\x_0$, for some small $\eps$. It follows that both the bias $\b(\x)$ and the bias gradient
\beq \label{eq:def B}
\B(\x) \triangleq \pd{\b}{\x}
\eeq
vanish at $\x = \x_0$. This formulation is the basis of the unconstrained unbiased CRB, a lower bound on the covariance at $\x_0$ which applies to all estimators whose bias gradient is zero at $\x_0$.

It turns out that even this requirement is too stringent in constrained settings. As we will see in Section~\ref{ss:sparse}, estimators of the coefficients of an overcomplete dictionary must have a nonzero bias gradient matrix. The reason is related to the fact that unbiasedness is required over the set $\Bxz$, which, in the overcomplete setting, has a higher dimension than the number of measurements.

However, it can be argued that one is not truly interested in the bias at all points in $\Bxz$, since many of these points violate the constraint set $\XX$. A reasonable compromise is to require unbiasedness over $\Bxz \cap \XX$, i.e., over the neighborhood of $\x_0$ restricted to the constraint set $\XX$. This leads to a weaker requirement on the bias gradient $\B$ at $\x_0$. Specifically, the derivatives of the bias need only be specified in directions which do not violate the constraints. The exact formulation of this requirement depends on the nature of the set $\XX$. In the following subsections, we will investigate various constraint sets and derive the corresponding requirements on the bias function.

It is worth emphasizing that the dependence of the CRB on the constraints is manifested through the class of estimators being considered, or more specifically, through the allowed estimators' bias gradient matrices. By contrast, the unconstrained CRB applies to estimators having a fully specified bias gradient matrix. Consequently, the constrained bound applies to a wider class of estimators, and is thus usually lower than the unconstrained version of the CRB\@. In other words, estimators which are unbiased in the constrained setting, and thus applicable to the unbiased constrained CRB, are likely to be biased in the unconstrained context. Since a wider class of estimators is considered by the constrained CRB, the resulting bound is lower, thus explaining the puzzling phenomenon described in the beginning of this subsection.

\subsection{Locally Balanced Constraints}
\label{ss:loc bal}

We now consider a class of constraint sets, called locally balanced sets, which encompass the sparsity constraints of Section~\ref{se:sparse backgnd}. Roughly speaking, a locally balanced set is one which is locally defined at each point by the directions along which one can move without leaving the set. Formally, a metric space $\XX$ is said to be locally balanced if, for all $\x \in \XX$, there exists an open set $\CC \subset \XX$ such that $\x \in \CC$ and such that, for all $\x' \in \CC$ and for all $|\lambda| \le 1$, we have
\beq \label{eq:def loc bal}
\x + \lambda (\x' - \x) \in \CC.
\eeq
As we will see, locally balanced sets are useful in the context of the constrained CRB, as they allow us to identify the feasible directions along which the bias gradient must be specified.

An example of a locally balanced set is given in Fig.~\ref{fi:loc bal}(a), which represents a union of two subspaces. In Fig.~\ref{fi:loc bal}(a), for any point $\x \in \XX$, and for any point $\x' \in \XX$ sufficiently close to $\x$, the entire line segment between $\x$ and $\x'$, as well as the line segment in the opposite direction, are also in $\XX$. This illustrates the fact that any union of subspaces is locally balanced, and, in particular, so are the sparse estimation settings of Section~\ref{se:sparse backgnd} \cite{eldar09, eldar09b, gedalyahu09}. As another example, consider any open set, such as the open ball in Fig.~\ref{fi:loc bal}(b). For such a set, any point $\x$ has a sufficiently small neighborhood $\CC$ such that, for any $\x' \in \CC$, the line segment connecting $\x$ to $\x'$ is contained in $\XX$. On the other hand, the curve in Fig.~\ref{fi:loc bal}(c) is not locally balanced, since the line connecting $\x$ to any other point on the set does not lie within the set.\footnote{We note in passing that since the curve in Fig.~\ref{fi:loc bal}(c) is continuously differentiable, it can be locally approximated by a locally balanced set. Our derivation of the CRB can be extended to such approximately locally balanced sets in a manner similar to that of \cite{gorman90}, but such an extension is not necessary for the purposes of this paper.}

Observe that the neighborhood of a point $\x$ in a locally balanced set $\XX$ is entirely determined by the set of feasible directions $\v$ along which infinitesimal changes of $\x$ do not violate the constraints. These are the directions $\v = \x' - \x$ for all points $\x' \ne \x$ in the set $\CC$ of \eqref{eq:def loc bal}. Recall that we seek a lower bound on the performance of estimators whose bias gradient is defined over the neighborhood of $\x_0$ restricted to the constraint set $\XX$. Suppose for concreteness that we are interested in unbiased estimators. For a locally balanced constraint set $\XX$, this implies that
\beq \label{eq:pre T-unbias}
\B \v = \zero
\eeq
for any feasible direction $\v$. In other words, all feasible directions must be in the nullspace of $\B$. This is a weaker condition than requiring the bias gradient to equal zero, and is thus more useful for constrained estimation problems. If an estimator $\hx$ satisfies \eqref{eq:pre T-unbias} for all feasible directions $\v$ at a certain point $\x_0$, we say that $\hx$ is $\XX$-unbiased at $\x_0$. This terminology emphasizes the fact that $\XX$-unbiasedness depends both on the point $\x_0$ and on the constraint set $\XX$.

Consider the subspace $\FF$ spanned by the feasible directions at a certain point $\x \in \XX$. We refer to $\FF$ as the feasible subspace at $\x$. Note that $\FF$ may include infeasible directions, if these are linear combinations of feasible directions. Nevertheless, because of the linearity of \eqref{eq:pre T-unbias}, any vector $\u \in \FF$ satisfies $\B\u = \zero$, even if $\u$ is infeasible. Thus, $\XX$-unbiasedness is actually a property of the feasible subspace $\FF$, rather than the set of feasible directions.

Since $\XX$ is a subset of a finite-dimensional Euclidean space, $\FF$ is also finite-dimensional, although different points in $\XX$ may yield subspaces having differing dimensions. Let $\u_1, \ldots, \u_l$ denote an orthonormal basis for $\FF$, and define the matrix
\beq \label{eq:def U}
\U = [ \u_1, \ldots, \u_l ].
\eeq
Note that $\u_i$ and $\U$ are functions of $\x$. For a given function $\x$, different orthonormal bases can be chosen, but the choice of a basis is arbitrary and will not affect our results.

As we have seen, $\XX$-unbiasedness at $\x_0$ can alternatively be written as $\B\u = \zero$ for all $\u \in \FF$, or, equivalently
\beq \label{eq:T-unbias}
\B\U = \zero.
\eeq
The constrained CRB can now be derived as a lower bound on all $\XX$-unbiased estimators, which is a weaker requirement than ``ordinary'' unbiasedness.

Just as $\XX$-unbiasedness was defined by requiring the bias gradient matrix to vanish when multiplied by any feasible direction vector, we can define $\XX$-biased estimators by requiring a specific value (not necessarily zero) for the bias gradient matrix when multiplied by a feasible direction vector. In an analogy to \eqref{eq:T-unbias}, this implies that one must define a value for the matrix $\B\U$. Our goal is thus to construct a lower bound on the covariance at a given $\x$ achievable by any estimator whose bias gradient $\B$ at $\x$ satisfies $\B\U = \P$, for a given matrix $\P$. This is referred to as specifying the $\XX$-bias of the estimator at $\x$.

\subsection{The CRB for Locally Balanced Constraints}

It is helpful at this point to compare our derivation with prior work on the constrained CRB, which considered continuously differentiable constraint sets of the form \eqref{eq:constr}. It has been previously shown \cite{gorman90} that inequality constraints of the type $\g(\x) \le \zero$ have no effect on the CRB\@. Consequently, we will consider constraints of the form
\beq \label{eq:eq constr}
\XX = \{ \x \in \RR^n: \f(\x)=\zero \}.
\eeq
Define the $k \times n$ matrix $\F(\x) = \partial \f / \partial \x$. For simplicity of notation, we will omit the dependence of $\F$ on $\x$. Assuming that the constraints are non-redundant, $\F$ is a full-rank matrix, and thus one can define an $n \times (n-k)$ matrix $\W$ (also dependent on $\x$) such that
\beq
\F\W = \zero, \quad \W^T\W = \I.
\eeq
The matrix $\W$ is closely related to the matrix $\U$ spanning the feasible direction subspace of locally balanced sets. Indeed, the column space $\Ra{\W}$ of $\W$ is the tangent space of $\XX$, i.e., the subspace of $\RR^n$ containing all vectors which are tangent to $\XX$ at the point $\x$. Thus, the vectors in $\Ra{\W}$ are precisely those directions along which infinitesimal motion from $\x$ does not violate the constraints, up to a first-order approximation. It follows that if a particular set $\XX$ is both locally balanced and continuously differentiable, its matrices $\U$ and $\W$ coincide. Note, however, that there exist sets which are locally balanced but not continuously differentiable (and vice versa).

With the above formulation, the CRB for continuously differentiable constraints can be stated as a function of the the matrix $\W$ and the bias gradient $\B$ \cite{ben-haim09}. In fact, the resulting bound depends on $\B$ only through $\B\W$. This is to be expected in light of the discussion of Section~\ref{ss:bias req}: The bias should be specified only for those directions which do not violate the constraint set. Furthermore, the proof of the CRB in \cite[Theorem~1]{ben-haim09} depends not on the formulation \eqref{eq:eq constr} of the constraint set, but merely on the class of bias functions under consideration. Consequently, one can state the bound without any reference to the underlying constraint set. To do so, let $\y$ be a measurement vector with pdf $p(\y;\x)$, which is assumed to be differentiable with respect to $\x$. The Fisher information matrix (FIM) $\J(\x)$ is defined as
\beq \label{eq:def J}
\J(\x) = \E{\bD \bD^T}
\eeq
where
\beq \label{eq:def bD}
\bD = \pd{\log p(\y;\x)}{\x}.
\eeq
We assume that the FIM is well-defined and finite. We further assume that integration with respect to $\y$ and differentiation with respect to $\x$ can be interchanged, a standard requirement for the CRB\@. We then have the following result.

\begin{theorem}
\label{th:crb}
Let $\hx$ be an estimator and let $\B = \partial \b / \partial \x$ denote the bias gradient matrix of $\hx$ at a given point $\x_0$. Let $\U$ be an orthonormal matrix, and suppose that $\B\U$ is known, but that $\B$ is otherwise arbitrary. If
\beq \label{eq:UUM in UUJUU}
\Ra{\U(\U+\B\U)^T)} \subseteq \Ra{\UUJUU}
\eeq
then the covariance of $\hx$ at $\x_0$ satisfies
\beq \label{eq:th:crb}
\Cov(\hx) \succeq (\U+\B\U) \left(\U^T\J\U \right)^\pinv (\U+\B\U)^T.
\eeq
Equality is achieved in \eqref{eq:th:crb} if and only if
\beq \label{eq:th:crb eq cond}
\hx = \x_0 + \b(\x_0) + (\U+\B\U) \left( \U^T\J\U \right)^\pinv \U^T \bD
\eeq
in the mean square sense, where $\bD$ is defined by \eqref{eq:def bD}. Conversely, if \eqref{eq:UUM in UUJUU} does not hold, then there exists no finite-variance estimator with the required bias gradient.
\end{theorem}

As required, no mention of constrained estimation is made in Theorem~\ref{th:crb}; instead, partial information about the bias gradient is assumed. Apart from this restatement, the theorem is identical to \cite[Theorem~1]{ben-haim09}, and its proof is unchanged. However, the above formulation is more general in that it can be applied to any constrained setting, once the constraints have been translated to bias gradient requirements. In particular, Theorem~\ref{th:crb} provides a CRB for locally balanced sets if the matrix $\U$ is chosen as a basis for the feasible direction subspace of Section~\ref{ss:loc bal}.

\section{Bounds on Sparse Estimation}
\label{se:sparse bounds}

In this section, we apply the CRB of Theorem~\ref{th:crb} to several sparse estimation scenarios. We begin with an analysis of the problem of estimating a sparse parameter vector.

\subsection{Estimating a Sparse Vector}
\label{ss:sparse}

Suppose we would like to estimate a parameter vector $\alf_0$, known to belong to the set $\TT$ of \eqref{eq:def T}, from measurements $\y$ given by \eqref{eq:y=H alf + w}. To determine the CRB in this setting, we begin by identifying the feasible subspaces $\FF$ corresponding to each of the elements in $\TT$. To this end, consider first vectors $\alf \in \TT$ for which $\|\alf\|_0 = s$, i.e., vectors having maximal support. Denote by $\{ i_1, \ldots, i_s \}$ the support set of $\alf$. Then, for all $\delta$, we have
\beq
\|\alf + \delta \e_{i_k}\|_0 = \|\alf\|_0 = s,
\quad  k=1,\ldots,s
\eeq
where $\e_j$ is the $j$th column of the identity matrix. Thus $\alf + \delta \e_{i_k} \in \TT$, and consequently, the vectors $\{ \e_{i_1}, \ldots, \e_{i_s} \}$ are all feasible directions, as is any linear combination of these vectors. On the other hand, for any $j \notin \supp(\alf)$ and for any nonzero $\delta$, we have $\|\alf + \delta \e_j\|_0 = s+1$, and thus $\e_j$ is not a feasible direction; neither is any other vector which is not in $\spn\{\e_{i_1}, \ldots, \e_{i_s}\}$. It follows that the feasible subspace $\FF$ for points having maximal support is given by $\spn\{\e_{i_1}, \ldots, \e_{i_s}\}$, and a possible choice for the matrix $\U$ of \eqref{eq:def U} is
\beq \label{eq:U when =s}
\U = [ \e_{i_1}, \ldots, \e_{i_s} ]
\quad \text{for } \|\alf\|_0 = s.
\eeq

The situation is different for points $\alf$ having $\|\alf\|_0 < s$. In this case, vectors $\e_i$ corresponding to \emph{any} direction $i$ are feasible directions, since
\beq
\|\alf + \delta \e_i\|_0 \le \|\alf\|_0 + 1 \le s.
\eeq
Because the feasible subspace is defined as the span of all feasible directions, we have
\beq
\FF \supseteq \spn\{ \e_1, \ldots, \e_p \} = \RR^p.
\eeq
It follows that $\FF = \RR^p$ and thus a convenient choice for the matrix $\U$ is
\beq \label{eq:U when <s}
\U = \I  \quad \text{for } \|\alf\|_0 < s.
\eeq
Consequently, whenever $\|\alf\|_0 < s$, a specification of the $\TT$-bias amounts to completely specifying the usual estimation bias $\b(\x)$.

To invoke Theorem~\ref{th:crb}, we must also determine the FIM $\J(\alf)$. Under our assumption of white Gaussian noise, $\J(\alf)$ is given by \cite[p.~85]{kay93}
\beq \label{eq:half J}
\J(\alf) = \frac{1}{\sigma^2} \H^T\H.
\eeq
Using \eqref{eq:U when =s}, \eqref{eq:U when <s}, and \eqref{eq:half J}, it is readily shown that
\beq \label{eq:half UJU}
\U^T\J\U =
\begin{cases}
\frac{1}{\sigma^2} \Hal^T \Hal & \text{when } \|\alf\|_0 = s \\
\frac{1}{\sigma^2} \H^T  \H    & \text{when } \|\alf\|_0 < s
\end{cases}
\eeq
where $\Hal$ is the $p \times s$ matrix consisting of the columns of $\H$ indexed by $\supp(\alf)$.

We now wish to determine under what conditions \eqref{eq:UUM in UUJUU} holds. Consider first points $\alf_0$ for which $\|\alf_0\|_0 = s$. Since, by \eqref{eq:spark req}, we have $\spark(\H)>s$, it follows that in this case $\U^T\J\U$ is invertible. Therefore
\beq
\Ra{\UUJUU} = \Ra{\U\U^T}.
\eeq
Since
\beq
\Ra{\U\U^T(\I+\B^T)} \subseteq \Ra{\U\U^T}
\eeq
we have that condition \eqref{eq:UUM in UUJUU} holds when $\|\alf_0\|_0=s$.

The condition \eqref{eq:UUM in UUJUU} is no longer guaranteed when $\|\alf_0\|_0 < s$. In this case, $\U=\I$, so that \eqref{eq:UUM in UUJUU} is equivalent to
\beq \label{eq:I+B in HH}
\Ra{\I+\B^T} \subseteq \Ra{\H^T\H}.
\eeq
Using the fact that $\Ra{\H^T\H} = \Ra{\H^T}$ and that, for any matrix $\Q$, $\Ra{\Q^T} = \Nu{\Q}^\perp$, we find that \eqref{eq:I+B in HH} is equivalent to
\beq \label{eq:N(H) in N(I+B)}
\Nu{\H} \subseteq \Nu{\I+\B}.
\eeq
Combining these conclusions with Theorem~\ref{th:crb} yields the following CRB for the problem of estimating a sparse vector.

\begin{theorem} \label{th:alf}
Consider the estimation problem \eqref{eq:y=H alf + w} with $\alf_0$ given by \eqref{eq:def T}, and assume that \eqref{eq:spark req} holds. For a finite-variance estimator $\half$ of $\alf_0$ to exist, its bias gradient matrix $\B$ must satisfy \eqref{eq:N(H) in N(I+B)} whenever $\|\alf_0\|_0 < s$. Furthermore, the covariance of any estimator whose $\TT$-bias gradient matrix is $\B\U$ satisfies
\begin{align} \label{eq:th:alf}
\Cov(\half) &\succeq \sigma^2 (\I+\B) (\H^T\H)^\pinv (\I+\B^T) \notag\\
  &\hspace{11em} \text{ when } \|\alf_0\|_0 < s, \notag\\
\Cov(\half) &\succeq \sigma^2 (\U+\B\U) (\Halz^T \Halz)^{-1} (\U+\B\U)^T \notag\\
  &\hspace{11em} \text{ when } \|\alf_0\|_0 = s.
\end{align}
Here, $\Halz$ is the matrix containing the columns of $\H$ corresponding to $\supp(\alf_0)$.
\end{theorem}

Let us examine Theorem~\ref{th:alf} separately in the underdetermined and well-determined cases. In the well-determined case, in which $\H$ has full row rank, the nullspace of $\H$ is trivial, so that \eqref{eq:N(H) in N(I+B)} always holds. It follows that the CRB is always finite, in the sense that we cannot rule out the existence of an estimator having any given bias function. Some insight can be obtained in this case by examining the $\TT$-unbiased case. Noting also that $\H^T\H$ is invertible in the well-determined case, the bound for $\TT$-unbiased estimators is given by
\begin{align} \label{eq:alf well-det unbiased}
\Cov(\half) &\succeq \sigma^2 (\H^T\H)^{-1}
  &\text{ when } \|\alf_0\|_0 &< s, \notag\\
\Cov(\half) &\succeq \sigma^2 \U (\Halz^T \Halz)^{-1} \U^T
  &\text{ when } \|\alf_0\|_0 &= s.
\end{align}

From this formulation, the behavior of the CRB can be described as follows. When $\alf_0$ has non-maximal support ($\|\alf_0\|_0 < s$), the CRB is identical to the bound which would have been obtained had there been no constraints in the problem. This is because $\U=\I$ in this case, so that $\TT$-unbiasedness and ordinary unbiasedness are equivalent. As we have seen in Section~\ref{ss:bias req}, the CRB is a function of the class of estimators under consideration, so the unconstrained and constrained bounds are equivalent in this situation. The bound $\sigma^2 (\H^T\H)^{-1}$ is achieved by the unconstrained LS estimator
\beq
\half = (\H^T\H)^{-1}\H^T\y
\eeq
which is the minimum variance unbiased estimator in the unconstrained case. Thus, we learn from Theorem~\ref{th:alf} that for values of $\alf_0$ having non-maximal support, no $\TT$-unbiased technique can outperform the standard LS estimator, which does not assume any knowledge about the constraint set $\TT$.

On the other hand, consider the case in which $\alf_0$ has maximal support, i.e., $\|\alf_0\|_0 = s$. Suppose first that $\supp(\alf_0)$ is known, so that one must estimate only the nonzero values of $\alf_0$. In this case, a reasonable approach is to use the oracle estimator \eqref{eq:def xo}, whose covariance matrix is given by $\sigma^2 \U (\Halz^T \Halz)^{-1} \U^T$ \cite{candes07}. Thus, when $\alf_0$ has maximal support, Theorem~\ref{th:alf} states that $\TT$-unbiased estimators can perform, at best, as well as the oracle estimator, which is equivalent to the LS approach when the support of $\alf_0$ is known.

The situation is similar, but somewhat more involved, in the underdetermined case. Here, the condition \eqref{eq:N(H) in N(I+B)} for the existence of an estimator having a given bias gradient matrix no longer automatically holds. To interpret this condition, it is helpful to introduce the mean gradient matrix $\M(\alf)$, defined as
\beq
\M(\alf) = \pd{\E{\half}}{\alf} = \I + \B.
\eeq
The matrix $\M(\alf)$ is a measure of the sensitivity of an estimator to changes in the parameter vector. For example, a $\TT$-unbiased estimator is sensitive to any \emph{feasible} change in $\alf$. Thus, $\Nu{\M}$ denotes the subspace of directions to which $\half$ is insensitive. Likewise, $\Nu{\H}$ is the subspace of directions for which a change in $\alf$ does not modify $\H\alf$. The condition \eqref{eq:N(H) in N(I+B)} therefore states that for an estimator to exist, it must be insensitive to changes in $\alf$ which are unobservable through $\H\alf$, at least when $\|\alf\|_0 < s$. No such requirement is imposed in the case $\|\alf\|_0 = s$, since in this case there are far fewer feasible directions.

The lower bound \eqref{eq:th:alf} is similarly a consequence of the wide range of feasible directions obtained when $\|\alf\|_0 < s$, as opposed to the tight constraints when $\|\alf\|_0 = s$. Specifically, when $\|\alf\|_0 < s$, a change to any component of $\alf$ is feasible and hence the lower bound equals that of an unconstrained estimation problem, with the FIM given by $\sigma^{-2} \H^T \H$. On the other hand, when $\|\alf\|_0 = s$, the bound is effectively that of an estimator with knowledge of the particular subspace to which $\alf$ belongs; for this subspace the FIM is the submatrix $\U^T\J\U$ given in \eqref{eq:half UJU}. This phenomenon is discussed further in Section~\ref{se:discuss}.

Another difference between the well-determined and underdetermined cases is that when $\H$ is underdetermined, an estimator cannot be $\TT$-unbiased for all $\alf$. To see this, recall from \eqref{eq:T-unbias} that $\TT$-unbiased estimators are defined by the fact that $\B\U=\zero$. When $\|\alf\|_0 < s$, we have $\U=\I$ and thus $\TT$-unbiasedness implies $\B=\zero$, so that $\Nu{\I+\B} = \{ \zero \}$. But since $\H$ is underdetermined, $\Nu{\H}$ is nontrivial. Consequently, \eqref{eq:N(H) in N(I+B)} cannot hold for $\TT$-unbiased estimators when $\|\alf\|_0 < s$.

The lack of $\TT$-unbiased estimators when $\|\alf_0\|_0 < s$ is a direct consequence of the fact that the feasible direction set at such $\alf_0$ contains all of the directions $\e_1, \ldots, \e_p$. The conclusion from Theorem~\ref{th:alf} is then that no estimator can be expected to be unbiased in such a high-dimensional neighborhood, just as unbiased estimation is impossible in the $p$-dimensional neighborhood $\Baz$, as explained in Section~\ref{ss:bias req}. However, it is still possible to obtain a finite CRB in this setting by further restricting the constraint set: if it is known that $\|\alf_0\|_0 = \tilde{s} < s$, then one can redefine $\TT$ in \eqref{eq:def T} by replacing $s$ with $\tilde{s}$. This will enlarge the class of estimators considered $\TT$-unbiased, and Theorem~\ref{th:alf} would then provide a finite lower bound on those estimators. Such estimators will not, however, be unbiased in the sense implied by the original constraint set.

While an estimator cannot be unbiased for \emph{all} $\alf \in \TT$, unbiasedness is possible at points $\alf$ for which $\|\alf\|_0 = s$. In this case, Theorem~\ref{th:alf} produces a bound on the MSE of a $\TT$-unbiased estimator, obtained by calculating the trace of \eqref{eq:th:alf} in the case $\B\U=\zero$. This bound is given by
\beq \label{eq:crb T-unbiased}
\E{\|\half - \alf_0\|_2^2} \ge \sigma^2 \Tr((\Halz^T \Halz)^{-1}), \quad
\|\alf_0\|_0 = s.
\eeq

The most striking feature of \eqref{eq:crb T-unbiased} is that it is identical to the oracle MSE \eqref{eq:oracle mse}. However, the CRB is of additional importance because of the fact that the ML estimator achieves the CRB in the limit when a large number of independent measurements are available, a situation which is equivalent in our setting to the limit $\sigma \rightarrow 0$. In other words, an MSE of \eqref{eq:crb T-unbiased} is achieved at high SNR by the ML approach \eqref{eq:ml}, as we will illustrate numerically in Section~\ref{se:numer}. While the ML approach is computationally intractable in the sparse estimation setting, it is still implementable in principle, as opposed to $\alfor$, which relies on unavailable information (namely, the support set of $\alf_0$). Thus, Theorem~\ref{th:crb} gives an alternative interpretation to comparisons of estimator performance with the oracle.

Observe that the bound \eqref{eq:crb T-unbiased} depends on the value of $\alf_0$ (through its support set, which defines $\Halz$). This implies that some values of $\alf_0$ are more difficult to estimate than others. For example, suppose the $\ell_2$ norms of some of the columns of $\H$ are significantly larger than the remaining columns. Measurements of a parameter $\alf_0$ whose support corresponds to the large-norm columns of $\H$ will then have a much higher SNR than measurements of a parameter corresponding to small-norm columns, and this will clearly affect the accuracy with which $\alf_0$ can be estimated. To analyze the behavior beyond this effect, it is common to consider the situation in which the columns $\h_i$ of $\H$ are normalized so that $\|\h_i\|_2 = 1$. In this case, for sufficiently incoherent dictionaries, $\Tr((\Halz^T \Halz)^{-1})$ is bounded above and below by a small constant times $s$, so that the CRB is similar for all values of $\alf_0$. To see this, let $\mu$ be the coherence of $\H$ \cite{tropp06}, defined (for $\H$ having normalized columns) as
\beq
\mu \triangleq \max_{i \ne j} \left| \h_i^T \h_j \right| .
\eeq
By the Gershgorin disc theorem, the eigenvalues of $\Halz^T \Halz$ are in the range $[1 - s\mu, 1 + s\mu]$. It follows that the unbiased CRB \eqref{eq:crb T-unbiased} is bounded above and below by
\beq
\frac{s\sigma^2}{1+s\mu} \le \sigma^2 \Tr((\Halz^T \Halz)^{-1}) \le \frac{s\sigma^2}{1-s\mu}.
\eeq
Thus, when $s$ is somewhat smaller than $1/\mu$, the CRB is roughly equal to $s \sigma^2$ for all values of $\alf_0$. As we have seen in Section~\ref{ss:est techniques}, for sufficiently small $s$, the worst-case MSE of practical estimators, such as BPDN and the DS, is $O(s \sigma^2 \log p)$. Thus, practical estimators come almost within a constant of the unbiased CRB, implying that they are close to optimal for all values of $\alf_0$, at least when compared with unbiased techniques.

\subsection{Denoising and Deblurring}
\label{ss:deblur}

We next consider the problem \eqref{eq:y=Ax+w}, in which it is required to estimate not the sparse vector $\alf_0$ itself, but rather the vector $\x_0 = \D \alf_0$, where $\D$ is a known dictionary matrix. Thus, $\x_0$ belongs to the set $\SS$ of \eqref{eq:def S}. We assume for concreteness that $\D$ has full row rank and that $\A$ has full column rank. This setting encompasses the denoising and deblurring problems described in Section~\ref{ss:sparse setting}, with the former arising when $\A=\I$ and the latter obtained when $\A$ represents a blurring kernel. Similar calculations can be carried out when $\A$ is rank-deficient, a situation which occurs, for example, in some interpolation problems.

Recall from Section~\ref{ss:sparse setting} the assumption that every $\x \in \SS$ has a \emph{unique} representation $\x = \D\alf$ for which $\alf$ is in the set $\TT$ of \eqref{eq:def T}. We denote by $\r(\cdot)$ the mapping from $\SS$ to $\TT$ which returns this representation. In other words, $\r(\x)$ is the unique vector in $\TT$ for which
\beq
\x = \D \r(\x) \quad \text{and} \quad \|\r(\x)\|_0 \le s.
\eeq
Note that while the mapping $\r$ is well-defined, actually calculating the value of $\r(\x)$ for a given vector $\x$ is, in general, NP-hard.

In the current setting, unlike the scenario of Section~\ref{ss:sparse}, it is always possible to construct an unbiased estimator. Indeed, even without imposing the constraint \eqref{eq:def S}, there exists an unbiased estimator. This is the LS or maximum likelihood estimator, given by
\beq
\hx = (\A^T\A)^{-1} \A^T \y.
\eeq
A standard calculation demonstrates that the covariance of $\hx$ is
\beq \label{eq:LS cov}
\sigma^2 (\A^T\A)^{-1}.
\eeq
On the other hand, the FIM for the setting \eqref{eq:y=Ax+w} is given by
\beq \label{eq:deb J}
\J = \frac{1}{\sigma^2} \A^T\A.
\eeq
Since $\A$ has full row rank, the FIM is invertible. Consequently, it is seen from \eqref{eq:LS cov} and \eqref{eq:deb J} that the LS approach achieves the CRB $\J^{-1}$ for unbiased estimators. This well-known property demonstrates that in the unconstrained setting, the LS technique is optimal among all unbiased estimators.

The LS estimator, like any unbiased approach, is also $\SS$-unbiased. However, with the addition of the constraint $\x_0 \in \SS$, one would expect to obtain improved performance. It is therefore of interest to obtain the CRB for the constrained setting. To this end, we first note that since $\J$ is invertible, we have $\Ra{\UUJUU} = \Ra{\U\U^T}$ for any $\U$, and consequently \eqref{eq:UUM in UUJUU} holds for any matrix $\B$. The bound \eqref{eq:th:crb} of Theorem~\ref{th:crb} thus applies regardless of the bias gradient matrix.

For simplicity, in the following we derive the CRB for $\SS$-unbiased estimators. A calculation for arbitrary $\SS$-bias functions can be performed along similar lines. Consider first values $\x \in \SS$ such that $\|\r(\x)\|_0 < s$. Then, $\|\r(\x) + \delta \e_i\|_0 \le s$ for any $\delta$ and for any $\e_i$. Therefore,
\beq
\x + \delta \D \e_i \in \SS
\eeq
for any $\delta$ and $\e_i$. In other words, the feasible directions include all columns of $\D$. Since it is assumed that $\D$ has full row rank, this implies that the feasible subspace $\FF$ equals $\RR^n$, and the matrix $\U$ of \eqref{eq:def U} can be chosen as $\U=\I$.

\begin{figure*}
\centerline{%
\subfigure[]{%
\includegraphics{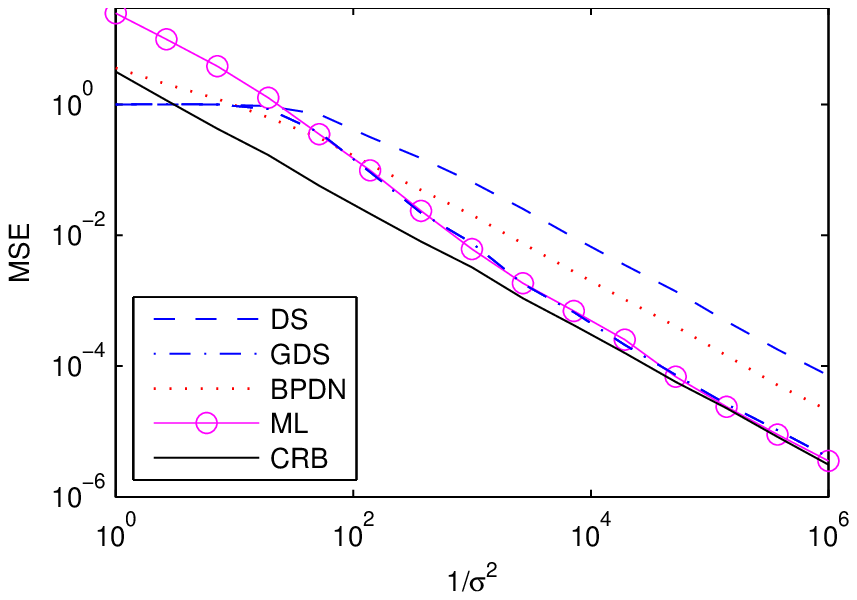}  
\label{fi:snr}}
\hfil %
\subfigure[]{%
\includegraphics{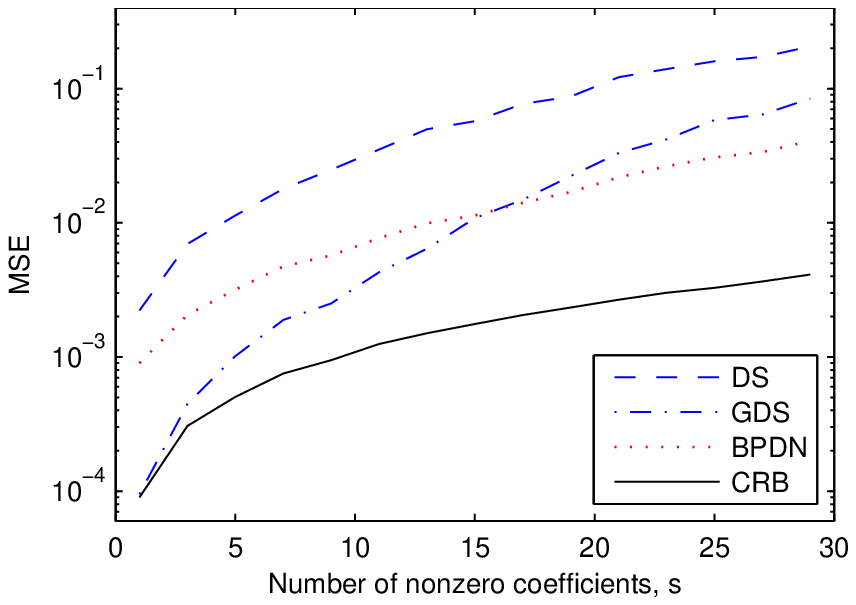}      
\label{fi:spar}}
}
\caption{MSE of various estimators compared with the unbiased CRB \eqref{eq:crb T-unbiased}, for (a) varying SNR and (b) varying sparsity levels.}
\label{fi:sim}
\end{figure*}

Next, consider values $\x \in \SS$ for which $\|\r(\x)\|_0 = s$. Then, for sufficiently small $\delta>0$, we have $\|\r(\x) + \delta \v\|_0 \le s$ if and only if $\v = \e_i$ for some $i \in \supp(\r(\x))$. Equivalently,
\beq
\x + \delta \v \in \SS \text{ if and only if } \v = \D\e_i \text{ and } i \in \supp(\r(\x)).
\eeq
Consequently, the feasible direction subspace in this case corresponds to the column space of the matrix $\D_\x$ containing the $s$ columns of $\D$ indexed by $\supp(\r(\x))$. From \eqref{eq:spark req D} we have $\spark(\D)>s$, and therefore the columns of $\D_\x$ are linearly independent. Thus the orthogonal projector onto $\FF$ is given by
\beq \label{eq:def P}
\P \triangleq \U\U^T = \D_\x (\D_\x^T \D_\x)^{-1} \D_\x^T.
\eeq
Combining these calculations with Theorem~\ref{th:crb} yields the following result.

\begin{theorem} \label{th:deblur}
Consider the estimation setting \eqref{eq:y=Ax+w} with the constraint \eqref{eq:def S}, and suppose $\spark(\D) > 2s$. Let $\hx$ be a finite-variance $\SS$-unbiased estimator. Then,
\begin{align}
\Cov(\hx) &\succeq \sigma^2 (\A^T \A)^{-1}
&\text{when } \|\r(\x)\|_0 < s, \notag\\
\Cov(\hx) &\succeq \sigma^2 \left( \P \A^T\A \P \right)^\pinv
&\text{when } \|\r(\x)\|_0 = s.
\label{eq:th:deblur}
\end{align}
Here, $\P$ is given by \eqref{eq:def P}, in which $\D_\x$ is the $n \times s$ matrix consisting of the columns of $\D$ participating in the (unique) $s$-element representation $\D\alf$ of $\x$.
\end{theorem}

As in Theorem~\ref{th:alf}, the bound exhibits a dichotomy between points having maximal and non-maximal support. In the former case, the CRB is equivalent to the bound obtained when the support set is known, whereas in the latter the bound is equivalent to an unconstrained CRB. This point is discussed further in Section~\ref{se:discuss}.

\section{Numerical Results}
\label{se:numer}

In this section, we demonstrate the use of the CRB for measuring the achievable MSE in the sparse estimation problem \eqref{eq:y=H alf + w}. To this end, a series of simulations was performed. In each simulation, a random $100 \times 200$ dictionary $\H$ was constructed from a zero-mean Gaussian IID distribution, whose columns $\h_i$ were normalized so that $\|\h_i\|_2=1$. A parameter $\alf_0$ was then selected by choosing a support uniformly at random and selecting the nonzero elements as Gaussian IID variables with mean $0$ and variance $1$. Noisy measurements $\y$ were obtained from \eqref{eq:y=H alf + w}, and $\alf_0$ was then estimated using BPDN \eqref{eq:bpdn}, the DS \eqref{eq:ds}, and the GDS \eqref{eq:gds}. The regularization parameters were chosen as $\tau = 2\sigma \sqrt{\log p}$ and $\gamma = 4\sigma \sqrt{\log(p-s)}$, rules of thumb which are motivated by a theoretical analysis \cite{ben-haim09c}. The MSE of each estimate was then calculated by repeating this process with different realizations of the random variables. The unbiased CRB was calculated using \eqref{eq:crb T-unbiased}. In this case, the unbiased CRB equals the MSE of the oracle estimator \eqref{eq:def xo}, but as we will see below, interpreting \eqref{eq:crb T-unbiased} as a bound on unbiased estimators provides further insight into the estimation problem.

A first set of experiments was conducted to examine the CRB at various SNR levels. In this simulation, the ML estimator \eqref{eq:ml} was also computed, in order to verify its convergence to the CRB at high SNR\@. Since the ML approach is computationally prohibitive when $p$ and $s$ are large, this necessitated the selection of the rather low support size $s=3$. The MSE and CRB were calculated for 15 SNR values by changing the noise standard deviation $\sigma$ between $1$ and $10^{-3}$. The MSE of the ML approach, as well as the other estimators of Section~\ref{ss:est techniques}, is compared with the CRB in Fig.~\ref{fi:snr}. The convergence of the ML estimator to the CRB is clearly visible in this figure. The performance of the GDS is also impressive, being as good or better than the ML approach. Apparently, at high SNR, the DS tends to correctly recover the true support set, in which case GDS \eqref{eq:gds} equals the oracle \eqref{eq:def xo}. Perhaps surprisingly, applying a LS estimate on the support set obtained by BPDN (which could be called a ``Gauss--BPDN'' strategy) does not work well at all, and in fact results in higher MSE than a direct application of BPDN. (The results for the Gauss--BPDN method are not plotted in Fig.~\ref{fi:sim}.)

Note that some estimation techniques outperform the oracle MSE (or CRB) at low SNR\@. It may appear surprising that a practical technique such as the DS outperforms the oracle. The explanation for this stems from the fact that the CRB \eqref{eq:crb T-unbiased} is a lower bound on the MSE of \emph{unbiased} estimators. The bias of most estimators tends to be negligible in low-noise settings, but often increases with the noise variance $\sigma^2$. Indeed, when $\sigma^2$ is as large as $\|\alf_0\|_2^2$, the measurements carry very little useful information about $\alf_0$, and an estimator can improve performance by shrinkage. Such a strategy, while clearly biased, yields lower MSE than a naive reliance on the noisy measurements. This is indeed the behavior of the DS and BPDN, since for large $\sigma^2$, the $\ell_1$ regularization becomes the dominant term, resulting in heavy shrinkage. Consequently, it is to be expected that such techniques will outperform even the best unbiased estimator at low SNR, as indeed occurs in Fig.~\ref{fi:snr}.

The performance of the estimators of Section~\ref{ss:est techniques}, excluding the ML method, was also compared for varying sparsity levels. To this end, the simulation was repeated for 15 support sizes in the range $1 \le s \le 30$, with a constant noise standard deviation of $\sigma = 0.01$. The results are plotted in Fig.~\ref{fi:spar}. While a substantial gap exists between the CRB and the MSE of the practical estimators in this case, the general trend in both cases describes a similar rate of increase as $s$ grows. Interestingly, a drawback of the GDS approach is visible in this setting: as $s$ increases, correct support recovery becomes more difficult, and shrinkage becomes a valuable asset for reducing the sensitivity of the estimate to random measurement fluctuations. The LS approach practiced by the GDS, which does not perform shrinkage, leads to gradual performance deterioration.

Results similar to Fig.~\ref{fi:sim} were obtained for a variety of related estimation scenarios, including several deterministic, rather than random, dictionaries $\H$.

\section{Discussion}
\label{se:discuss}

In this paper, we extended the CRB to constraint sets satisfying the local balance condition (Theorem~\ref{th:crb}). This enabled us to derive lower bounds on the achievable performance in various estimation problems (Theorems \ref{th:alf} and~\ref{th:deblur}). In simple terms, Theorems \ref{th:alf} and~\ref{th:deblur} can be summarized as follows. The behavior of the CRB differs depending on whether or not the parameter has maximal support (i.e., $\|\alf\|_0 = s$). In the case of maximal support, the bound equals that which would be obtained if the sparsity pattern were known; this can be considered an ``oracle bound''. On the other hand, when $\|\alf\|_0 < s$, performance is identical to the unconstrained case, and the bound is substantially higher. We now discuss some practical implications of these conclusions. To simplify the discussion, we consider the case of unbiased estimators, though analogous conclusions can be drawn for any bias function.

When $\|\alf\|_0 = s$ and all nonzero elements of $\alf$ are considerably larger than the standard deviation of the noise, the support set can be recovered correctly with high probability (at least if computational considerations are ignored). Thus, in this case an estimator can mimic the behavior of the oracle, and the CRB is expected to be tight. Indeed, in the high SNR limit, the ML estimator achieves the unbiased CRB\@. On the other hand, when the support of $\alf$ is not maximal, the unbiasedness requirement demands sensitivity to changes in all components of $\alf$, and consequently the bound coincides with the unconstrained CRB\@. Thus, as claimed in Section~\ref{se:crb}, in underdetermined cases no estimator is unbiased for all $\alf \in \SS$.

An interesting observation can also be made concerning maximal-support points $\alf$ for which some of the nonzero elements are close to zero. The CRB in this ``low-SNR'' case corresponds to the oracle MSE, but as we will see, the bound is loose for such values of $\alf$. Intuitively, at low-SNR points, any attempt to recover the sparsity pattern will occasionally fail. Consequently, despite the optimistic CRB, it is unlikely that the oracle MSE can be achieved. Indeed, the covariance matrix of any finite-variance estimator is a continuous function of $\alf$ \cite{lehmann98}, and the fact that performance is bounded by the (much higher) unconstrained bound when $\|\alf\|_0 < s$ implies that performance must be similarly poor for low SNR\@.

This excessive optimism is a result of the local nature of the CRB\@: The bound is a function of the estimation setting only in an $\eps$-neighborhood of the parameter itself. Indeed, the CRB depends on the constraint set only through the feasible directions, which were defined in Section~\ref{ss:loc bal} as those directions which do not violate the constraints for \emph{sufficiently small} deviations. Thus, for the CRB, it is entirely irrelevant if some of the components of $\alf$ are close to zero, as long as $\supp(\alf)$ is held constant.

A tighter bound for sparse estimation problems may be obtained using the Hammersley--Chapman--Robbins (HCR) approach \cite{Hammersley50, ChapmanRobbins51, gorman90}, which depends on the constraints at points beyond the local neighborhood of $\x$. Such a bound is likely to yield tighter results for low SNR values, and will create a smooth transition between the regions of maximal and non-maximal support. However, the bound will depend on more complex properties of the estimation setting, such as the distance between $\D\alf$ and feasible points with differing supports. The derivation of such a bound is a subject for further research.

\section*{Acknowledgement}

The authors would like to thank Yaniv Plan for helpful discussions.
The authors are also grateful to the anonymous reviewers for their comments, which considerably improved the presentation of the paper.

\bibliographystyle{IEEEtran}
\bibliography{IEEEabrv,zvika}

\end{document}